\documentclass{elsart}
\usepackage{ifpdf}
\usepackage{graphicx,natbib,amssymb,lineno}
\usepackage {amsmath, latexsym, amscd, amssymb}
\usepackage {graphicx}
\usepackage[latin1]{inputenc}
\usepackage [french]{babel}
\ifpdf
\usepackage[%
  pdftitle={Groupes lin\'eaires finis permutant deux fois transitivement un ensemble  de droites},%
  pdfauthor={Lucas Vienne},%
  pdfsubject={The preprint document class elsart},%
  pdfkeywords={automorphisme,  groupe, graphe, graphe fortement r\'egulier, groupe doublement transitif, graphe de Paley},%
  pdfstartview=FitH,%
  bookmarks=true,%
  bookmarksopen=true,%
  breaklinks=true,%
  colorlinks=true,%
  linkcolor=blue,anchorcolor=blue,%
  citecolor=blue,filecolor=blue,%
  menucolor=blue,pagecolor=blue,%
  urlcolor=blue]{hyperref}
\else
\usepackage[%
  breaklinks=true,%
  colorlinks=true,%
  linkcolor=blue,anchorcolor=blue,%
  citecolor=blue,filecolor=blue,%
  menucolor=blue,pagecolor=blue,%
  urlcolor=blue]{hyperref}
\fi

\makeatletter
\def\elsartstyle{%
    \def\normalsize{\@setfontsize\normalsize\@xiipt{14.5}}
    \def\small{\@setfontsize\small\@xipt{13.6}}
    \let\footnotesize=\small
    \def\large{\@setfontsize\large\@xivpt{18}}
    \def\Large{\@setfontsize\Large\@xviipt{22}}
    \skip\@mpfootins = 18\p@ \@plus 2\p@
    \normalsize
}
\@ifundefined{square}{}{\let\Box\square}
\makeatother

\newcommand{\PPP}{\mathbb{P}}

\newcommand{\C}{\mathbb{C}} 
\newcommand{\R}{\mathbb{R}}

\newcommand{\F}{\mathbb{F}}

\newcommand{\EE}{\mathcal{E}}

\newcommand{\GG}{\mathcal{G}}

\newcommand{\JJ}{\mathcal{J}}

\newcommand{\PP}{\mathcal{P}}

\newcommand{\ov}{\overline}

\newcommand{\x}{\times}
\newcommand{\wt}{\widetilde}

\newcommand{\al}{\langle}
\newcommand{\ar}{\rangle}

\newcommand{\rmsi}{\mathrm{si}}

\newcommand{\Gram}{\mathrm{Gram}}

\newcommand{\nlb}{\nolinebreak}

\newcommand{\mm}{m\^eme }
\newcommand{\elt}{\'el\'ement }

\newcommand{\si }{si et seulement si }

\newcommand{\Dem }{{\it D\'emonstration. }}

\newcommand{\ttt}{th\'eor\`eme }

\newcommand{\oo}{\circ }

\newcommand{\ii }{\"{\i}}

\newcommand{\ev }{espace vectoriel }

\newcommand{\seps }{sous-espaces propres }

\newcommand{\ps }{ produit scalaire }

\newcommand{\bon }{base orthonorm\'ee }

\newcommand{\ga }{\alpha }

\newcommand{\gc }{\gamma }
\newcommand{\gd }{\delta }

\newcommand{\gve }{\varepsilon }
\newcommand{\gvf }{\varphi }

\newcommand{\gl}{\lambda}

\newcommand{\go }{\omega }
\newcommand{\gs}{\sigma}

\newcommand{\gD }{\Delta }

\newtheorem{lemf}{Lemme}
\newtheorem{thmf}{Th\'eor\`eme }

\newtheorem{propf}{Proposition}

\begin{document}

\begin{frontmatter}
\title{Groupes lin\'eaires finis permutant deux fois transitivement un ensemble de droites}

\author{Lucas Vienne}
\address{Departement de math\'ematiques. Universit\'e d'Angers. France}

\ead{lucas.vienne@univ-angers.fr}

\begin{abstract}
Let $n\geq 2$ be a positive integer, and $G$ a doubly transitive subgroup of the symmetric group on $X=\{1,\ldots,n\}$. In this paper we find all linear group representations $\rho$ of $G$ on an euclidean vector space $V$ which contains a set of equiangular vector lines $\GG=\{\al v_1\ar,\ldots,\al v_n\ar\}$ such that :\\
\hspace*{0,5cm}$(1)$ $V$ is generated by $v_1,\ldots,v_n$,\\
\hspace*{0,5cm}$(2)$ for all $i$ in $X$ and all $\gs$ in $G$, \quad $\al \rho_\gs(v_i)\ar=\al v_{\gs(i)}\ar$.\\[0,1cm]
Then we illustrate our construction when $G=SL_d(q)$, $q$ odd and $d\geq 2$.
 \end{abstract}

\end{frontmatter}
\section{Introduction} 
Un entier $n\geq 2$\  est donn\'e ainsi qu'un groupe de permutations $G$  sur l'ensemble \mbox{$X=\{1,\ldots,n\}$}. Nous dirons (par abus)  qu'une repr\'esentation lin\'eaire $\gc$ de $G$ dans un espace vectoriel $V$ est deux fois transitive  s'il existe une gerbe de droites vectorielles $\GG =\{\al v_1\ar,\ldots, \al v_n\ar\}$ de $V$ telle que :\\
(1) \ L'espace  $V$ est engendr\'e par  les vecteurs $v_1,\ldots,v_n.$\\
(2) \ Pour tout indice $i$ dans $X$ et tout $\gs$ dans $G$ on a $\al \gc_\gs(v_i)\ar=\al v_{\gs(i)}\ar$.

{\it  Notations}\\
1. Lorsque le contexte est sans ambigu\ii t\'e sur l'action $\gc$ d'un groupe $G$ sur un ensemble $Y$, on utilise, pour $\gs$ dans $G$ et $y$ dans $Y$, la notation simplifi\'ee $\gc_\gs(y)= \gs(y)$. Par exemple,  la relation (2) pr\'ec\'edente s'\'ecrit plus simplement \\[0,1cm]
\vspace{0,2cm}\centerline{ $\gs\al v_i\ar =\al \gs(v_i)\ar=\al v_{\gs(i)}\ar$.}
2. Le groupe $G$ \'etant fini, son image $\gc(G) $ est contenue dans
le groupe orthogonal  d'un \ps  sur $V$
qui sera not\'e $\gvf=(.|.)$  dans la suite.

{\it Quelques remarques}\\[0,1cm]
 Choisissons des g\'en\'erateurs $v_1,\ldots,v_n$ de \mm norme pour les droites de la gerbe $\GG$. 
La double transitivit\'e de $G$  nous montre qu'il existe des scalaires $\go, c$ (o\`u $\go\neq 0$) et des coefficients $\gve_{i,j}$ ($i\neq j $, dans $X$), tous dans $\{-1,+1\}$ tels que \\[0,1cm]
$(R1)$ \vspace{0,2cm}\centerline{$\forall i,j \in X, \quad (v_i | v_j)=
 \left\lbrace\begin{array}{ll}
 \go & \ \rmsi \ i=j \\
 \gve_{i,j}.c & \  \rmsi \ i\neq j
 \end{array}\right.$.}
De plus, pour chaque \elt $\gs$ dans $G$, il existe une liste de coefficients  $\nu^\gs=(\nu^\gs_1,\ldots, \nu^\gs_n)$ tous pris dans $\{-1,+1\}$ tels que \\[0,2cm]
$(R2)$ \vspace{0,2cm}\centerline{$\forall j \in X,\quad  \gs(v_j)=\nu^\gs_j.v_{\gs(j)}$.} 
Utilisant la relation $(R2)$, il vient pour $i\neq j$, \\[0,3cm]
 \vspace{0,3cm}\centerline{$\gve_{i,j}.c=(v_i|v_j)=(\gs(v_i)|\gs(v_j))=(\nu^\gs_i.v_{\gs(i)}|\nu^\gs_j.v_{\gs(j)})=\nu^\gs_i.\nu^\gs_j.\gve_{\gs(i),\gs(j)}.c$,}
  $(R3)$ \  donc,  si $c\neq 0$,  \vspace{0,3cm}\centerline{$\gve_{\gs(i),\gs(j)}=\nu^\gs_i.\nu^\gs_j. \gve_{i,j}.$\hspace{4cm}}
La double transitivit\'e de $G$ nous montre 
  que les coefficients $\gve_{i,j}$  sont d\'etermin\'es par la
  connaissance de l'un d'eux et  des coefficients $\nu_j^\gs$
  introduits en  $(R2)$. Notant $\EE=(\gve_{i,j})$  la matrice dont la diagonale principale est nulle $(\gve_{i,i}=0$) et les
   autres coefficients sont les $\gve_{i,j}$, d\'esignant par $I_n$ , la matrice de Gram du syst\`eme de vecteurs  $(v_1,\ldots, v_n)$ est
 donc de la forme  \\[0,2cm]
$(R4)$ \vspace{0,2cm}\centerline{ $\Gram(v_1,\ldots,
   v_n)=S_\EE(\go,c)=\go.I_n+c.\EE$.}
  
Dans cet article nous d\'eterminons quels choix des syst\`emes de coefficients  \\[0,1cm]
(3)\vspace{0,1cm}\centerline{$\nu^\gs=(\nu^\gs_1,\ldots, \nu^\gs_n)$,}
  satisfaisant aux  relations $(R2)$ d\'efinissent une repr\'esentation lin\'eaire de $G$ sur $V$, puis nous \'etudions  la nature, irr\'eductible ou non, de ces repr\'esentations.
   Pour \'enoncer nos r\'esultats,  introduisons le groupe multiplicatif produit \\[0,1cm]
\vspace{0,1cm}\centerline{$\Pi^n=\{\nu=(\nu_1,\ldots,\nu_n) \  |  \  \forall j,\ \nu_j=\pm 1\}$,}
 sur lequel on fait agir  $G$ ; pour tout $\nu=(\nu_1,\ldots,\nu_n)$ dans $\Pi^n$ et toute permutation $\gs$ de $X$ notons $\nu_\gs=(\nu_{\gs(1)},\ldots, \nu_{\gs(n)})$, puis posons \\[0,2cm]
$(4)$ \vspace{0,2cm}\centerline{$\forall (\gs,\nu) \in G\x \Pi^n, \quad \theta_\gs(\nu)=\nu_{\gs^{-1}}$.}
La relation imm\'ediate $(\nu.\nu')_\gs=\nu_\gs.\nu'_\gs$ nous montre que $\theta_\gs$ est un morphisme. Noter que dans la d\'efinition $(4)$, l'exposant $(-1)$ est l\`a pour s'assurer que $\theta$ un morphisme et non un anti-morphisme ($\theta_{\gs.\gs'}=\theta_{\gs}.\theta_{\gs'}$).
Notons $G\ltimes \Pi^n$ le produit semi-direct de $G$ par $\Pi^n$ sous
cette action.

Les th\'eor\`emes 1 et 2 donnent nos principaux r\'esultats. Nous les ferons suivre par quelques exemples.

\begin{thmf} $\ $\\
Soit $G$ un groupe de permutation $2$-transitif sur l'ensemble $X=\{1,\ldots,n\}$, 
$(e_1,\ldots,e_n)$ une base orthonorm\'ee d'un \ev  $E$ et \mbox{$\JJ : G\to G\ltimes \Pi^n$}  un morphisme de groupes de la forme \\[0,1cm]
$(5)$\vspace{0,2cm}\centerline{$\forall \gs \in G, \quad \JJ(\gs)=(\gs,\nu^\gs)$}
$1.$ Il existe une unique repr\'esentation lin\'eaire $\rho :G\to GL(E)$ dite {\rm associ\'ee} au morphisme $\JJ$, satisfaisant aux relations $(R2)$ :\\[0,1cm]
\vspace{0,2cm}\centerline{$\forall j \in X,\quad  \rho_\gs(e_j)=\nu^\gs_j.e_{\gs(j)}$.} 
$2.$ Le nombre de matrices \ $\EE=(\gve_{i,j})$,  sym\'etriques, \`a coefficients dans $\{\pm 1\}$
 et de diagonale principale nulle satisfaisant aux relations $(R3)$ \\[0,2cm]
\vspace{0,2cm}\centerline{$\forall \gs\in G,\forall i,j \in X, \quad \gve_{\gs(i),\gs(j)}=\nu^\gs_i.\nu^\gs_j. \gve_{i,j}$}
est \  $0$ ou $2$. Si $\EE$ est l'une d'elles l'autre est $-\EE$.\\
$3.$ Si la repr\'esentation $\rho :G\to GL(E)$ n'est pas
irr\'eductible, alors :\\ 
\hspace*{0,5cm}$a.$ Il existe une matrice $\EE$ satisfaisant aux conditions $(R3)$.\\
\hspace*{0,5cm}$b.$ Cette matrice poss\`ede exactement deux valeurs propres $\gl_1$
et $\gl_2$.\\
\hspace*{0,5cm}$c.$  Les \seps $E_1$ et $E_2$ associ\'es sont des sous-$G$-modules irr\'eductibles
 de $E$ sur lesquels le groupe $G$
induit des repr\'esentations doublement
transitives.
\end{thmf}

\begin{thmf} $\ $\\
Soit $G$ un groupe de permutation de  l'ensemble $X=\{1,\ldots,n\}$, et une repr\'esentation lin\'eaire $\gc :
G\to GL(V)$   doublement transitive de
$(G,X)$ sur un \ev euclidien $V$ de dimension finie.\\
Il existe un morphisme \mbox{$\JJ : G\to G\ltimes \Pi^n$} de la forme $(5)$ tel que :\\
$1.$ Si $\dim(V)=n$, la repr\'esentation $\gc$ est semblable \`a la repr\'esentation $\rho$ associ\'ee au morphisme $\JJ$.\\
$2.$  Si $\dim(V)<n$,  la repr\'esentation $\gc$ est semblable \`a l'une des deux
composantes  irr\'eductibles de  $\rho$ d\'ecrites dans la partie $3$ du \ttt $1$.

 \end{thmf}

Nous appliquerons ensuite ces th\'eor\`emes aux sous groupes $G$ de $GL_d(q)$ contenant $SL_d(q)$, agissant deux fois transitivement sur l'espace projectif $\PPP_{d-1}(q)$.

\section{Preuves} 

{\bf 2.1 \  R\'esultats  pr\'eliminaires} 

\begin{propf}$\ $\\ 
Soit $E$ un espace euclidien r\'eel muni d'une \bon $(e)=(e_1,\ldots,e_n)$, $S$ une matrice sym\'etrique r\'eelle de type $n\x n$, et pour tout couple de scalaires $(\go,c)$, posons $S(\go,c)=\go.I_n+c.S$. Notons   $\psi_{k,c}$ et $\wt{\psi}_{k,c}$ la forme bilin\'eaire et l'endomorphisme sym\'etrique associ\'es   \`a la matrice  $S(\go,c)$ dans la base $(e)$.\\
 On suppose $c\neq 0$.\\
 $1.a.$ L'application affine $\ga : \gl \to \go+c.\gl$ \'etablit une bijection entre les valeurs propres de $S=S(0,1)$ et celles de $S(\go,c)$.\\
$1.b.$ Les sous-espaces propres de $\wt{\psi}_{k,c}$ sont ind\'ependants du couple $(\go,c)$. Plus pr\'ecis\'ement, si  $E_{k,c}(\mu)$ d\'esigne le sous-espace propre de  $\wt{\psi}_{\go,c}$ associ\'e \`a la valeur propre $\mu$ on a : \\[-0,2cm]
\vspace{0,2cm}\centerline{ $E_{\go,c}(\go+c.\gl)=E_{0,1}(\gl)$}
$2.$ Soit $E=E_1\oplus \ldots\oplus E_s$ la d\'ecomposition de $E$ en somme des sous-espaces propres de $S$ associ\'es aux valeurs propres $\gl_1,\ldots,\gl_s$ qu'on suppose rang\'ees en fonction de leurs multiplicit\'es croissantes $m_1\leq \ldots \leq m_s$. 
Le noyau de $\wt{\psi}_{\go,c}$ est $(0)$ ou l'un des $E_i$
$(1\leq i\leq s)$. \\
$3.$ Le rang de $\wt{\psi}_{\go,c}$ est toujours sup\'erieur ou \'egal \`a $n-m_s$.
\end{propf}

\Dem Simple, on la laisse au lecteur. \hfill $\Box$

On se place dans la situation d\'ecrite dans
l'introduction : le groupe $G$ op\`ere deux fois transitivement sur
l'ensemble $X=\{1,\ldots,n\}$. La repr\'esentation lin\'eaire  $\gc
: G\to GL(V)$ est doublement transitive dans un \ev euclidien
$V$ et il existe des g\'en\'erateurs $v_1,\ldots,v_n$ de $V$ et des
scalaires $\nu^\gs_i$ (pour $(\gs,i)\in G\x X$) tels que les relations
$(R2)$ soient satisfaites. Le lemme suivant nous informe sur les
sous-modules de $V$.

\begin{lemf}$\ $\\
 Supposons que $U$ est un sous-$G$-module  de $V$ distinct de $V$.\\
 $1.$ Il existe un  projecteur $p_U$ de $V$ sur $U$ qui commute aux
 op\'erations de $G$ :\\[0,1cm]
\vspace{0,1cm}\centerline{$\forall \gs \in G,\quad p\oo \gs=\gs\oo
  p.$}
$2.$ Notant, pour tout indice $i\in X$, $u_i=p(v_i)$, alors \\[0,1cm]
\vspace{0,1cm}\centerline{$\forall \gs \in G,\  \gs(u_i)=\nu^\gs_iu_{\gs(i)}$.}
$3.$ Il existe une matrice $\EE=(\gve_{i,j})$ de diagonale principale
nulle, dont les  coefficients, pour $i\neq j$, sont dans $\{-1,+1\}$  et satisfont  aux relations $(R3)$ \\[0,1cm]
\vspace{0,2cm}\centerline{$\forall \gs\in G,\forall i,j \in X, \quad \gve_{\gs(i),\gs(j)}=\nu^\gs_i.\nu^\gs_j. \gve_{i,j}$}
$4.$ Pour des constantes non nulles convenablement choisies $\go_U$ et
$c_U$ on a \\[0,1cm]
\vspace{0,1cm}\centerline{
  $\Gram(u_1,\ldots,u_n)=S_\EE(\go_U,c_U)=\go_U I_n+c_U\EE$}
\end{lemf}

\Dem  La premi\`ere affirmation vient d'un r\'esultat classique : tout
sous-$G$-module $U$ de $V$ admet un suppl\'ementaire $W$ qui est aussi
un $G$-module. Alors les projecteurs $p_U$ et $p_W$ associ\'es \`a la
d\'ecomposition $V=U\oplus W$ commutent aux op\'erations de $G$. 
La deuxi\`eme affirmation en d\'ecoule directement puisque pour tout
indice $i\in X$, l'\'egalit\'e $\gs(v_i)=\nu^\gs_iv_{\gs(i)}$ conduit \`a\\[0,1cm]
\vspace{0,1cm}\centerline{ $\gs(u_i)=\gs\oo p(v_i)=p\oo
\gs(v_i)=p(\nu^\gs_iv_{\gs(i)})=\nu^\gs_ip(v_{\gs(i)})=\nu^\gs_iu_{\gs(i)}.$}
On en d\'eduit que les vecteurs $u_i$ sont tous de \mm norme et de la double transitivit\'e de $G$ sur $X$, 
 on d\'eduit
l'existence des constantes $\go_U$ et $c_U$ telles que le syst\`eme de
vecteurs $(u_1,\ldots,u_n)$ satisfasse aux relations $(R1)$. Mais le
coefficient $\go_U$ est non nul puisque c'est le carr\'e des vecteurs
non nuls $u_i$, et le coefficient $c_U$ est aussi non nul puisque le
rang de la matrice $\Gram(u_1,\ldots,u_n)$ \'etant inf\'erieur
\`a la dimension du module $U$, il est strictement inf\'erieur \`a
$n$.
Finalement la matrice $\EE=(\gve_{i,j})$ satisfait aux relations $(R3)$, si bien que  $\Gram_\psi(u_1,\ldots,u_n)=S_\EE(\go_U,c_U)$
comme voulu. \hfill $\Box$

{\bf 2.2 \  Preuve du \ttt 1} %

On adopte les hypoth\`eses et notations de ce \ttt.\\
1. Pour chaque \elt $\gs$ de $G$ la relation $(R2)$ :\\[0,1cm]
\vspace{0,2cm}\centerline{$\forall j \in X,\quad  \rho_\gs(e_j)=\nu^\gs_j.e_{\gs(j)}$} 
d\'efinit une fonction $ \rho_\gs$ envoyant la base orthonorm\'ee
$(e)=(e_1,\ldots,e_n)$ sur une \bon de $E$, qui s'\'etend donc en une
isom\'etrie lin\'eaire de $E$. Ceci montre bien s\^ur l'unicit\'e de
$\rho$, mais il nous  reste \`a v\'erifier que $\rho$ est un morphisme.
Or, comme $\JJ$ est un morphisme on a \\[0,2cm]
\vspace{0,2cm}\centerline{ 
$\JJ(\gs. \gd)=\gs \gd.\nu^{\gs \gd}=\JJ(\gs).\JJ( \gd)=\gs\nu^\gs.\gd\nu^\gd=\gs\gd.\nu^\gs_\gd\nu^\gd$,}
donc \vspace{0,1cm}\centerline{ $\nu^{\gs \gd}=\nu^\gs_\gd\nu^\gd.$ \hspace{2cm}}
Par ailleurs : \  $\rho_{\gs \gd}(e_i)=\nu^{\gs.\gd}_ie_{\gs. \gd(i)}$
\  et \  
$\rho_{\gs}. \rho_{\gd}(e_i)=\rho_{\gs}(\nu^\gd_i.e_{\gd(i)})=\nu^\gd_i. \nu^\gs_{\gd(i)}e_{\gs. \gd(i)},$\\[0,2cm]
D'o\`u l'\'egalit\'e recherch\'ee, \  $\rho_{\gs\gd}(e_i)=\rho_{\gs}. \rho_{\gd}(e_i)$.

2. Le groupe $G$ op\'erant deux fois transitivement sur l'ensemble
$X$, il est clair, d'apr\`es la relation $(R3)$, que la connaissance
d'un des coefficients $\gve_{i,j}$ d\'etermine tous les autres. Par
ailleurs si $\EE$ est une matrice satisfaisant aux relations $(R3)$ il
en est bien s\^ur de \mm de $-\EE$, ce qui prouve l'assertion 2.

$3. a.$  Pour tout sous-$G$-module $U$ de $E$, distinct de $E$, 
 le lemme $1$  nous montre qu'il existe une matrice $\EE$ satisfaisant aux relations $(R3)$, ce qui prouve $a.$
Dans la suite, si $(u_1,\dots,u_n)$  est l'image  de la base orthonorm\'ee $(e_1,\dots,e_n)$ par la projection orthogonale $p_U : E\to U$  sa matrice de Gram sera simplement not\'ee $\Gram_U=\Gram(u_1,\dots,u_n)$. \\
Toujours d'apr\`es le lemme 1, la matrice  $\Gram_U$  est de la forme\\[0,1cm]
\vspace{0,1cm}\centerline{ $\Gram_U=S(\go,c)=\go
I_n+c.\EE$}
 pour des constantes non nulles $\go $ et $c$.
 L'endomorphisme sym\'etrique $\wt{\psi}_{\go,c}$ de matrice
 $S(\go,c)$ dans la base $(e)=(e_1,\ldots,e_n)$ se diagonalise sur
 $E$ et ses \seps sont clairement des sous-$G$-modules de $E$, or
 d'apr\`es la proposition $1$, ce sont les m\^emes que ceux de la
 matrice $\EE=S(0,1)$. Soient $\gl_1,\ldots,\gl_s$ les valeurs propres
 de  $\EE$, rang\'ees dans l'ordre de leur multiplicit\'es croissantes
 $m_1\leq \cdots \leq m_s$, et $E_1,,\ldots,E_s$ les sous-espaces
 propres associ\'es. 
  Remarquons que $\EE$ n'\'etant pas une
 homoth\'etie, on a $s\geq 2$. Comme les $E_u$ sont des sous-$G$-modules de
 $E$, le projecteur orthogonal $p_u$ de $E$ sur $E_u$ ($1\leq u\leq s$)
 commute avec l'action de $G$ donc, toujours d'apr\`es le lemme 1,  la matrice $\Gram_{E_u}$ est  de la forme $\Gram(u_1,\ldots,u_n)=S_\EE(\go_u,c_u)$ pour des constantes
 non nulles $\go_u$ et $c_u$. D'apr\`es la proposition 1, le noyau de
 $S_\EE(\go_u,c_u)$ est l'un des espaces propres $E_t$ ($1\leq t\leq s$) de $\EE$ et son 
rang est sup\'erieur ou \'egal \`a $n-m_s$. Or ce rang est aussi
inf\'erieur ou \'egal \`a la dimension $m_u$ de $E_u$. Donc  \\[0,1cm]
\vspace{0,2cm}\centerline{  $n-m_s\leq m_u$ \  soit  \ $n\leq m_u +m_s$.}
Mais comme $s\geq 2$ on peut choisir $u\neq s$ et il vient alors \\[0,1cm]
\vspace{0,2cm}\centerline{   $n\leq m_u +m_s \leq n$,}
donc $n=m_u +m_s$, ce qui implique $s=2, u=1$, et prouve $b$.

Le module $E$ est donc la somme directe orthogonale des \seps $E_1$ et $E_2$ de la matrice $\EE$  de dimensions respectives $m_1$ et $m_2$ (o\`u $m_1\leq m_2$) qui sont des sous-$G$-modules de $E$.  Pour achever la d\'emonstration du \ttt 1, il nous faut montrer que les sous-modules $E_1$ et $E_2$ sont irr\'eductibles. \\
Soit $W$ un sous-module irr\'eductible de $E$. D'apr\`es le lemme 1, sa matrice de Gram  est de la forme \\[0,1cm]
\vspace{0,2cm}\centerline{$\Gram_W= S(\go_W,c_W)=\go_WI_n+c_W\EE$,}
et, n'\'etant pas de rang $n$, elle est li\'ee \`a $G_{E_1}$ ou $G_{E_2}$.  Son rang $m=\dim(W)$ est  donc \'egal \`a $m_1$ ou $m_2$, ce qui prouve d\'ej\`a que $E_1$ est irr\'eductible. Raisonnant par l'absurde,  si $E_2$ ne l'est pas, il se d\'ecompose en somme orthogonale de sous-modules $E_2=E_2^1\oplus\cdots \oplus E_2^k$ qui sont tous de dimension $m_1$ et admettent  des matrices de Gram  de rang $m_1$, donc li\'ees \`a la matrice  $\Gram_{E_1}$ de $E_1$. L'espace $E$ \'etant somme orthogonale de $E_1$ et des $E_2^j$ ($1\leq j\leq k$), il est facile de v\'erifier que sa matrice de Gram, qui n'est autre que l'identit\'e, est la somme des matrices de Gram de chacun des termes de la somme directe. Or cette somme  bien s\^ur une matrice li\'ee \`a $\Gram_{E_1}$ qui ne peut pas \^etre l'identit\'e. Ceci prouve, par l'absurde, que $E_2$ n'est pas r\'eductible et ach\`eve  la d\'emonstration du \ttt 1.\hfill $\Box$

{\bf 2.3 \  Preuve du \ttt 2} 

D'apr\`es l'introduction, on peut choisir  des vecteurs $v_1,\ldots,v_n$ dans $V$ tels que :\\
\hspace*{0,3cm}* $V=\al v_1,\ldots, v_n\ar$.\\
\hspace*{0,3cm}* Le groupe $G$ op\`ere $2$-transitivement la gerbe de droites $\GG=\{\al v_1\ar,\ldots, \al v_n\ar\}$.\\
\hspace*{0,3cm}* Les relations $(R2)$, $(R3)$ et $(R4)$  sont satisfaites, pour des constantes $\go,c$ et une matrice $\EE=(\gve_{i,j})$ convenablement choisies.\\
Reprenant un raisonnement fait plus haut on voit que pour tout indice $i\in X$,\\[0,1cm]
\vspace{0,21cm}\centerline{$\gc_{\gs \gd}(v_i)=\nu^{\gs.\gd}_iv_{\gs. \gd(i)}$ \quad et \quad 
$\gc_{\gs}. \gc_{\gd}(v_i)=\gc_{\gs}(\nu^\gd_i.v_{\gd(i)})=\nu^\gd_i. \nu^\gs_{\gd(i)}v_{\gs. \gd(i)},$}
d'o\`u l'on d\'eduit que  $\nu^{\gs \gd}=\nu^\gs_\gd\nu^\gd$, puis que l'application $\JJ : \gs \to \gs.\nu^\gs$ est un morphisme de $G$ dans le produit semi-direct $G\ltimes \Pi^n$ puisque \\[0,1cm]
\vspace{0,2cm}\centerline{ 
$\JJ(\gs. \gd)=\gs \gd.\nu^{\gs \gd}$ \quad  et \quad  $\JJ(\gs).\JJ( \gd)=\gs\nu^\gs.\gd\nu^\gd=\gs\gd.\nu^\gs_\gd\nu^\gd$.}
Soit  $\rho : G\to GL(E)$ la repr\'esentation  associ\'ee au morphisme
$\JJ$ (cf. thm 1). L'espace  $E$ est donc euclidien, muni d'une \bon
$e_1,\ldots,e_n$ telle que  la repr\'esentation $\rho$ satisfasse aux relations $(R2)$ :\\[0,1cm]
\vspace{0,2cm}\centerline{$\forall j \in X,\quad  \rho_\gs(e_j)=\nu^\gs_j.e_{\rho_\gs(j)}$.} 
Comme la repr\'esentation $\gc$ satisfait aux m\^emes  relations, l'application \mbox{$\pi : e_j \to v_j$} (pour $j\in X$) se prolonge en un morphisme du $G$-module $E$ sur le $G$-module \nlb$V$. \\
1. Si $\dim V=n$, \ $\pi $ est un isomorphisme.\\
2. Sinon le noyau de $\pi$ est un sous-module de $E$ qui ne peut qu'\^etre  $E$ ou l'un des deux sous-modules $E_1$ et $E_2$ de $E$ d\'ecrits dans le \ttt 1. Son image est donc aussi soit nulle, soit isomorphe \`a l'un de ces deux modules. \hfill $\Box$

\section{Quelques exemples} 

{\bf \  Exemple 1} 

Voici un exemple de repr\'esentation lin\'eaire $\gc
: G \to GL(E)$, associ\'ee \`a un morphisme $\JJ : G\to G\ltimes \Pi^n$,
qui est doublement transitive et  irr\'eductible. \\
Nous utiliserons le r\'esultat classique suivant (voir  par exemple $[1]$) :

\begin{lemf} $\ $\\
Soit $G$ un groupe de permutation de  $ Y=\{1,\ldots,m\}$ et $\rho$ la
repr\'esentation lin\'eaire de $G$ sur $\R^m$  associ\'ee qui, pour une \bon $(e_1,\ldots,e_m)$ s'\'ecrit  : \\
\vspace{0,2cm}\centerline{$\forall \gs\in G, \forall i\in Y,\quad \rho_\gs(e_i)=e_{\gs(i)}$.}
Notons  $\chi$ le caract\`ere de la repr\'esentation $\rho$,  $\chi_1,\ldots,\chi_k, $  \ les caract\`eres irr\'eductibles du groupe $G$, et
pour pour chaque indice $j$ $(1\leq j \leq k)$, soit  $\mu_j=(\chi | \chi_j)$ la multiplicit\'e
du caract\`ere $\chi_j$ dans la repr\'esentation $\rho$.\\
Alors la somme  $\sum_j\mu_j^2$ est \'egale au nombre d'orbites de $G$ dans son action naturelle
 sur $Y\x Y$.
\end{lemf} 

Appliquons ce r\'esultat  dans la situation suivante : \\
Soit  $\F_3$ le corps \`a trois \'el\'ements, $V$ l'\ev $\F_3^3$,  et $G$ le groupe $SL(V)$ agissant naturellement sur l'ensemble $Y=V^*$ des $26$ vecteurs non nuls de $V$.\\
Il est facile de v\'erifier que le groupe $G$ poss\`ede trois orbites dans son action naturelle sur $Y\x Y$ qui sont :\\
* L'ensemble des couples de vecteurs  $(v,v)$ pour $v\in V^*$.\\
* L'ensemble des couples de vecteurs $(v,-v)$ pour $v\in V^*$.\\
* L'ensemble des couples de vecteurs  $(v,w)$ pour $v,w\in V^*,$ non colin\'eaires.\\
Choisissons un syst\`eme de repr\'esentants $v_1,\ldots v_{13}$ des droites vectorielles de $V$ et notons $\ov{v}_1,\ldots \ov{v}_{13}$ leurs oppos\'es. La repr\'esentation lin\'eaire $\gc : G\to GL(\R^{26})$ associ\'ee \`a l'action de naturelle de $G$ sur $Y=V^*=\{v_1,\ldots v_{13},\ov{v}_1,\ldots \ov{v}_{13}\}$ permute la base $Y$ de $\R^{26}$.
D'apr\`es le lemme 2, cette repr\'esentation $\gc $ se d\'ecompose en une somme de repr\'esentations irr\'eductibles $\gc_j$ intervenant avec des multiplicit\'es $\mu_j$ telles que $3=\sum_j \mu_j^2$. On en d\'eduit imm\'ediatement que $\gc$ est une somme de trois repr\'esentations irr\'eductibles non semblables. Or il n'est pas difficile  de les d\'ecrire  car la repr\'esentation $\gc$ admet clairement comme sous-repr\'esentation l'espace $U$ engendr\'e par les vecteurs de la forme $v_i+\ov{v}_i$, dans lequel la droite vectorielle $\sum_i v_i+\sum_i \ov{v}_i$ $(1\leq i\leq 13)$ est elle \mm un sous espace stable par $G$. Il s'ensuit que l'orthogonal $E$ de $U$ dans $V$ est lui-\mm un sous espace stable par $G$ et irr\'eductible, puisque $V$ est somme d'exactement trois sous $G$-modules irr\'eductibles. Regardons comment $G$ agit sur $E$.
Comme $E=U^\perp$, il admet pour base l'ensemble des 13 vecteurs de la forme $e_i=v_i-\ov{v}_i$, et l'action de $G$ sur $E$ est donn\'ee par des relations \\[0.1cm]
\vspace{0,2cm}\centerline{$\forall \gs \in G,  \quad \gs(e_i)=\gs(v_i-\ov{v}_i)=\nu^\gs_i(v_{\gs(i)}-\ov{v}_{\gs(i)})= \nu^\gs_i.e_{\gs(i)}$,}
o\`u les coefficients $\nu^\gs_i$ valent tous $\pm 1$.  
Donc la repr\'esentation  de $G=SL_3(3)$ sur $E$ donn\'ee par ces relations  est irr\'eductible et associ\'ee au morphisme \mbox{$\JJ :G \to G\ltimes \Pi^n$} donn\'e par $\JJ(\gs)=\gs.\nu^\gs$.

{\bf \ Exemple 2} %
 
G\'en\'eralisons maintenant l'id\'ee utilis\'ee dans l'exemple 1.

Soit $\F_q$ un  corps fini \`a $q$ \'el\'ements, de caract\'eristique impaire, $C$
 l'ensemble de ses carr\'es non nuls et $\ov{C}$ l'ensemble
de ses non-carr\'es.  Choisissons un entier $d\geq 2$, notons $V$
l'\ev $F_q^d$, $n$ le nombre de droites vectorielles de $V$,  chacune
d'elle \'etant engendr\'ee par un vecteur $v_i$ $(i\in X=\{1,\dots,n\})$.
Soit enfin $G$ un sous-groupe de $GL(V)$ contenant le groupe sp\'ecial
lin\'eaire $SL(V)$. L'action naturelle de $SL(V)$ sur l'espace projectif $\PPP(V)$ \'etant doublement transitive, il en de \mm de celle de $G$. Elle se transmet aussi \`a $X$ en notant $\gs(i)$, pour tout indice  $i$ dans $X$ et tout $\gs$ dans $G$, l'unique indice  tel que \\[0.1cm]
 \vspace{0,2cm}\centerline{$\gs(\al v_i\ar)=\al v_{\gs(i)}\ar.$}
Plus pr\'ecis\'ement, il existe  un scalaire non nul $\gl_i$ tel que \\[0.1cm]
 \vspace{0,2cm}\centerline{$\gs(v_i)=\gl_i.v_{\gs(i)}.$}
  On d\'efinit alors une repr\'esentation $\gc$ de $G$ comme groupe de permutations de  l'ensemble $Y=\{v_1,-v_1,\ldots, v_n, -v_n\}$ en posant \\[0.1cm]
(6) \hspace*{0,5cm}$\begin{array}{lll}
* \ \gc_\gs(v_i)=\nu^\gs_i.v_{\gs(i)} \quad & \mathrm{ou} &
  \quad  \nu^\gs_i=\left\lbrace\begin{array}{rl}
1\  & \rmsi \quad \gl_i \in C\\
-1\  & \rmsi \quad  \gl_i \in \ov{C}
\end{array}\right.,\\
*\ \gc_\gs(-v_i)=-\gc_\gs(v_i).&&
\end{array}$

 Pour \'eviter toute confusion posons  maintenant pour chaque indice $i$, $\ov{v_i}=-v_i$, notons encore $\gc$ la repr\'esentation lin\'eaire du groupe $G$ sur l'espace $\R^{2n}$ associ\'ee \`a l'action du groupe $G$ sur $Y$ et soit $\chi$ son caract\`ere. D'apr\`es le lemme 2, le nombre $\go$ d'orbites de $G$ dans son action naturelle sur $Y\x Y$ vaut \\[0,1cm]
 \vspace{0,2cm}\centerline{$\displaystyle{\go= \sum_j\mu_j^2}$,  \quad   o\`u  \quad  $\mu_j=(\chi | \chi_j)$} 
 repr\'esente la multiplicit\'e du caract\`ere irr\'eductible $\chi_j$   dans la repr\'esentation \nlb$\gc$. 
 Regardons comment $G$ agit sur $\R^{2n}$. \\
 L'\'egalit\'e $\{\gc_\gs(v_i),\gc_\gs(\ov{v}_i)\}=\{v_{\gs(i)}, \ov{v}_{\gs(i)} \}$ valable pour tout $\gs$ dans $G$ nous montre que le
 sous-espace $U$ engendr\'e par les vecteurs  $u_i=v_i+\ov{v}_i$ \ 
 $(i\in X)$ est stable par l'action de $G$, ainsi que la droite vectorielle $D$ engendr\'ee par le vecteur  $\sum_iu_i,$ et son orthogonal dans $U$, $D^\perp\cap U$. L'orthogonal $E$ de $U$ dans $\R^{2n}$ est lui-\mm stable par $G$, engendr\'e par les $n$ vecteurs $e_i=v_i-\ov{v}_i$, et l'action de $G$ sur $E$ est donn\'ee par les relations \\[0.1cm]
(7) \vspace{0,2cm}\centerline{$\forall \gs \in G,  \quad \gc_\gs(e_i)=\gc_\gs(v_i-\ov{v}_i)=\nu^\gs_i(v_{\gs(i)}-\ov{v}_{\gs(i)})= \nu^\gs_i.e_{\gs(i)}$}
Donc la repr\'esentation  de $G$ sur $E$ donn\'ee par ces relations  est  associ\'ee au morphisme $\JJ :G \to G\ltimes \Pi^n$ donn\'e par $\JJ(\gs)=\gs.\nu^\gs$. 
Le  \ttt suivant d\'ecrit, suivant les diff\'erentes situations, la d\'ecomposition de $\gc$ en somme de repr\'esentations irr\'eductibles.

\begin{thmf}$\ $\\ 
Sous les hypoth\`eses pr\'ec\'edentes, soit $\gc  : G\to GL(E)$ la repr\'esentation de $G$ donn\'ee par les relations $(7)$ \\[0.1cm]
 \vspace{0,1cm}\centerline{$\forall i\in X, \forall \gs \in G,  \quad \gc_\gs(e_i)=\nu^\gs_i.e_{\gs(i)}$}
$-$ Si $d\geq 3$ la repr\'esentation $\gc$ est irr\'eductible.\\
$-$ Supposons $d=2$ et notons $GL_2^+(q)$ le sous-groupe de $GL_2(q)$
form\'e des automorphismes dont le d\'eterminant est un carr\'e dans $F_q$.\\
 \hspace*{0,3cm} $(a)$ Si $G\not\subset GL_2^+(q)$, la repr\'esentation
$\gc$ est irr\'eductible. \\ 
 \hspace*{0,3cm} $(b)$  Si $G \subset GL_2^+(q)$,  la repr\'esentation $\gc$ est r\'eductible.\\
Dans le cas $(b)$ le  $G$-module $E$ se d\'ecompose  en somme directe $E=E_1\oplus E_2$ de deux
sous-modules irr\'eductibles de \mm dimension, \'echang\'es par
l'action sur $E$ du groupe $GL_2(q)$. De plus, \\
$*$ si $q \not\equiv 1 \mod
4$, le module $E$ reste irr\'eductible sur $\R$, \\
$*$ si $q\equiv 1 \mod
4$, la d\'ecomposition \`a lieu sur $\R$ et les repr\'esentations de
$G$ sur $E_1$ et $E_2$ sont associ\'ees aux graphes de Paley $\PP(q)$
$($voir $[5])$.\\
\end{thmf}

\Dem  
Sous les hypoth\`eses g\'en\'erales du \ttt 3,  l'action de $G$ sur $X$ \'etant doublement transitive,  le \mbox{\ttt 1} nous dit que la repr\'esentation $\gc$ de $G$ sur $E$ est soit irr\'eductible sur $\R$, soit  somme directe de deux repr\'esentations irr\'eductibles, non semblables. Mais d'apr\`es le lemme 2, dans le premier cas le groupe $G$ poss\`ede trois orbites sur $Y\x Y$ tandis que dans le second il en poss\`ede quatre. Or l'action de $G$ sur $Y\x Y$ poss\`ede au moins trois parties stables qui sont :\\
* L'ensemble $\gD$ des couples de vecteurs  $(v,v)$ pour $v\in Y$.\\
* L'ensemble $\gD'$  des couples de vecteurs $(v,-v)$ pour $v\in Y$.\\
* L'ensemble $\nabla $ des couples de vecteurs  $(v,w)$ pour $v,w\in Y,$ non colin\'eaires.\\
De plus l'action $\gc$ de $G$ sur $Y$ est n\'ecessairement transitive sinon $\gD$ et $\gD'$ seraient form\'ees d'au moins deux orbites,  si bien que l'action $\gc\x \gc$ de $G$ sur  $Y\x Y$ aurait au moins cinq orbites, ce qui est impossible. Donc seule la 
partie $\nabla $ peut \^etre form\'ee d'une ou deux orbites sous l'action $\gc\x \gc$ de $G$. \\
Discutons suivant les valeurs de la dimension $d$ de l'espace $V$.

$-$ {\it Supposons  $d\geq 3$.}

 L'action naturelle du groupe $SL_d(q)$ sur l'ensemble des couples $(v,w)$ de vecteurs non li\'es dans $V$  est transitive, tout comme celle de $G$, donc l'action $\gc\x \gc$ de $G$ sur $\nabla $ est transitive. Le lemme 2 nous montre que dans ce cas la repr\'esentation $\gc$ de $G$ sur $E$ est irr\'eductible.
  
$-$ {\it Supposons  $d=2$.} 

\begin{lemf}$\ $\\
Le nombre $\theta$ d'orbites du groupe $G$, agissant sur l'ensemble $\nabla $ des couples de vecteurs $(v,w)$ pour $v,w$ dans $ Y$ non colin\'eaires est \\
 \hspace*{0,3cm} $*$ \ $\theta=2,$ \  si \ $G\subset GL_2^+(q)$,\\
 \hspace*{0,3cm} $*$ \ $\theta=1,$  \ si \ $G\not\subset GL_2^+(q)$.\\
\end{lemf}

\Dem  On sait d\'ej\`a que $G$ poss\`ede une ou deux orbites dans son
action $\gc\x \gc$ sur $\nabla$. De plus comme $G$ op\`ere transitivement sur les couples de droites $(\al v_i \ar, \al v_j\ar)$ distinctes, il  op\`ere transitivement sur $\nabla$ \si l'orbite de $(v_i,v_j)$ sous l'action de $G$ par $\gc\x\gc$ contient $(v_i,-v_j)$, $(-v_i,v_j)$ et $(-v_i,-v_j)$, ce qui revient \`a dire que $G$ contient un \elt $\gs$ tel que $\gs(v_i)=v_i$ et $\gs(v_j)=\gl.v_j$ o\`u le scalaire $\gl$ n'est pas un carr\'e (d'apr\`es la d\'efinition $(5)$ de l'action $\gc$). Mais comme $G$ contient le groupe $SL_2(q)$, ceci signifie aussi que $G$ contient un \elt dont le d\'eterminant n'est pas un carr\'e, ce qui prouve le lemme. \hfill $\Box$

Pour achever la preuve du \ttt 3, on se place maintenant dans le cas $(b)$ de ce th\'eor\`eme pour lequel $d=2$ et  $G\subset GL_2^+(q)$. \\
D'apr\`es le lemme 2, l'action $\gc$  de $G$ sur $E$ est r\'eductible, et plus pr\'ecis\'ement $E$ est somme de deux sous-$G$-modules irr\'eductibles $E_1$ et $E_2$. 
Mais, toujours, d'apr\`es le lemme 2,  l'action $\gc\x \gc$ du groupe lin\'eaire $GL_2(q)$ \'etant transitive sur l'ensemble $\nabla, $  la repr\'esentation de $GL_2(q)$ sur $E$ est irr\'eductible. Or  le sous-groupe $G$ est invariant dans $GL_2(q)$ puisqu'il contient $SL_2(q)$, donc  $GL_2(q)$ permute les deux sous-espaces $E_1$ et $E_2$. Distinguons  deux cas :

{\it Premier cas : $q\not\equiv 1 \mod 4$.}\\
L'automorphisme $\gs$ qui envoie la base $(v_1,v_2)$ de $V$ sur $(v_2,-v_1)$ est dans $SL_2(q)$  donc dans $G$, et comme $-1$ n'est pas un carr\'e dans $\F_q$ les relations $(6)$ et $(7)$ nous conduisent \`a $\gc_\gs(v_1)=v_2$, $\gc_\gs(v_2)=-v_1$ puis \\[0,1cm]
\vspace{0,2cm}\centerline{$\gc_\gs(e_1)=e_2=\nu^\gs_1e_2$ \quad et \quad $\gc_\gs(e_2)=-e_1=\nu^\gs_2$.}
S'il existait une matrice $\EE=(\gve_{i,j})$ satisfaisant aux relations $(R4)$  du \mbox{th\'eor\`eme 1}, on devrait avoir  \\[0,1cm]
\vspace{0,2cm}\centerline{$ \gve_{\gs(1),\gs(2)}=\gve_{2,1}=\gve_{1,2}=\nu^\gs_1\nu^\gs_2.\gve_{1,2}=-\gve_{1,2}$,}
ce qui est faux.  Donc $\EE$ n'existe pas et d'apr\`es le \ttt 1 la
repr\'esentation  $\gc$ de $G$ sur $E$ est irr\'eductible sur $\R$.
Sa d\'ecomposition  en somme directe $E=E_1\oplus E_2$ a donc lieu sur
le corps $\C$ des complexes.

{\it Deuxi\`eme cas : $q\equiv 1 \mod 4$.}\\
La d\'ecomposition $E=E_1\oplus E_2$ a maintenant lieu sur le corps des r\'eels et les repr\'esentations de $G$ sur $E_1$ et $E_2$ sont associ\'ees \`a des familles \'equiangulaires de droites li\'ees aux graphes de Paley $\PP(q)$ (voir $[5]$ pour une construction de ces familles de droites). \hfill $\Box$

\section{Bibliographie}

\label{thebib}

{\it Articles}\\

[1] \quad Dominique de Caen, {\it Large equiangular sets of lines in Euclidean spaces.} The electronic journal of comninatorics. November 9, 2000.

[2] \quad L. Nguyen Van The, {\it On a problem of Specker about Euclidean
representations of finite graphs.}, (2008)\\
http://arxiv.org/abs/0810.2359

[3] \quad Aidan Roy , {\it Minimal Euclidean representations of graphs}. (2008)\\
 http://arxiv.org/abs/0812.3707

[4] \quad L. Vienne, {\it Repr\'esentations lin\'eaires des graphes finis.} (2009)\\ 
	http://arxiv.org/abs/0902.1874.
	
[5] \quad L. Vienne, {\it Groupes d'isom\'etries permutant doublement transitivement un ensemble de droites vectorielles (1).} (2009)\\ 
 http://fr.arxiv.org/abs/0903.0912.

[6] \quad  Brouwer, A. E.; Cohen, A. M.; and Neumaier, A. {\it Conference Matrices and Paley Graphs.} In Distance Regular Graphs. New York: Springer-Verlag.

{\it Livres}

$[7]$ \quad Godsil, Chris, and Royle, {\it  Algebraic Graph Theory.} \mbox{New York: Springer. (2001)} 

$[8]$ \quad G\'erard Rauch, {\it Les groupes finis et leurs repr\'esentations.} \mbox{Ellipses. (2001)}

$[9]$\quad  Norman Biggs. {\it Finite groups of automorphisms.} \mbox{London Mathematical Society. (1970)}

$[10]$\quad  Norman Biggs. {\it Algebraic Graph Theory.} \mbox{Cambridge University Press,  (1993)}

$[11]$\quad P. J. Cameron and J.H.Van Lint. {\it Designs, graphs,codes and their links, vol 22 of London Mathematical Society Student Texts,} Cambridge University Press, Cambridge, (1991).

\end{document}